\documentclass[12pt]{amsart}
\usepackage{amssymb}
\usepackage{verbatim}
\usepackage{xspace}
\usepackage{pstricks,pst-node}
\usepackage{vmargin}
\setpapersize{USletter}
\setmarginsrb{1.1in}{1.2in}{1.1in}{1.2in}{12 pt}{25 pt}{12 pt}{30 pt}

\theoremstyle{plain}
\newtheorem*{theorem}{Theorem}

\newtheorem*{proposition}{Proposition}

\theoremstyle{definition}

\theoremstyle{remark}

\newcommand{\ds}{\displaystyle}
\newcommand{\cstar}{\ensuremath{\text{C}^{*}}-}

\newcommand{\ssstar}[2]{s^{\phantom{*}}_{#1}s^*_{#2}}
\newcommand{\sstars}[2]{s^*_{#1}s^{\phantom{*}}_{#2}}

\DeclareMathOperator{\Obj}{Obj}
\newcommand{\bbN}{\mathbb N}
\newcommand{\bbZ}{\mathbb Z}
\newcommand{\bbT}{\mathbb T}
\newcommand{\sA}{\mathcal A}
\newcommand{\sB}{\mathcal B}
\newcommand{\sC}{\mathcal C}
\newcommand{\sD}{\mathcal D}
\newcommand{\sG}{\mathcal G}

\begin{document}
\title[Spectral Theorem for Bimodules]{The Spectral Theorem for
  Bimodules in\\
   Higher Rank Graph \cstar algebras}
\author{Alan Hopenwasser}
\address{Department of Mathematics\\
        University of Alabama\\
        Tuscaloosa, AL 35487}
\email{ahopenwa@euler.math.ua.edu}

\date{\today}

 \begin{abstract}
In this note we extend the spectral theorem for bimodules to the
higher rank graph {\ensuremath{\text{C}^{*}}-}algebra context.  Under
the assumption that the graph is row finite and has no sources, we
show that a bimodule over a natural abelian subalgebra is determined
by its spectrum iff it is generated by the Cuntz-Krieger partial
isometries which it contains iff the bimodule is invariant under the
gauge automorphisms.  We also show that the natural abelian subalgebra
is a masa iff the higher rank graph satisfies an aperiodicity condition.
 \end{abstract}
\maketitle

\section{Introduction}\label{S:intro}

Many \cstar algebras can be coordinitized -- a property that proves
very useful both in the study of the \cstar algebra and also of its
subalgebras. Coordinitization is achieved by presenting the \cstar
algebra as a groupoid \cstar algebra.  The unit space of the groupoid
is associated with an abelian subalgebra which is often, though not
always, a masa.  (The abelian subalgebra depends on the choice of
coordinates and need not be intrinsic.)  A great many of the
(non-self-adjoint) subalgebras of a groupoid \cstar algebra either
contain the ``diagonal'' abelian algebra or are a bimodule over it.
When the groupoid is $r$-discrete and principal, one of the most
fundamental tools used in the study of subalgebras is the spectral
theorem for bimodules of Muhly and Solel \cite{MR90m:46098}.  Roughly
speaking, this says that a bimodule is determined by the coordinates
on which it is supported.

When the groupoid is not principal, it is no longer true that a
bimodule is determined by its spectrum.  For graph \cstar algebras,
\cite{hpp} contains a characterization of those bimodules which are
determined by their spectra: these are the bimodules which are
invariant under the gauge automorphisms.  (Another equivalent
condition is that the bimodule be generated by the Cuntz-Krieger
partial isometries which it contains.)  Graph \cstar algebras have
been extensively studied in the last decade; see \cite{raeburn} for an 
excellent summary and a bibliography of relevant papers.  More
recently, considerable attention has turned to a multi-dimensional
analog, the higher rank graph \cstar algebras.

In the paper in which higher rank graph \cstar algebras were first
formalized \cite{MR2001b:46102}, Kumjian and Pask modified the path
groupoid model for graph \cstar algebras to produce a model for higher 
rank graph \cstar algebras.  The purpose of this note is to extend the 
spectral theorem for bimodules as it appears in \cite{hpp} for graph
\cstar algebras to the higher rank context.  Section \ref{S:grp} will
provide a brief review of the notation and construction of higher rank 
graph \cstar algebras and their associated path groupoids.  Section
\ref{S:stb} is devoted to the spectral theorem for bimodules in the
higher rank context.  It also contains a characterization of when the
``diagonal'' is a masa.

\section{Higher rank \cstar
  algebras and the path groupoid}\label{S:grp}

A $k$-graph $(\Lambda,d)$ is a small category $\Lambda$ together with
a functor $d \colon \Lambda \to \bbN^k$ which satisfies the following
factorization property: if $\lambda \in \Lambda$ and 
$d(\lambda) = m+n$ with $m,n \in \bbN^k$, then there exist unique
$\mu, \nu \in \Lambda$ such that $\lambda = \mu \nu$, $d(\mu) =m$, and 
$d(\nu)= n$.  For $n \in \bbN^k$, we let 
$\Lambda^n = d^{-1}(n)$ and note that $\Lambda^0$ can be identified
with the objects in $\Lambda$.

When $k=1$, $\Lambda$ is the category of finite paths from a directed
graph; $\Lambda^0$ is the set of vertices; $\Lambda^1$ is the set of
directed edges; and $\Lambda^n$ is the set of paths of length $n$.  A
higher rank graph is a multi-dimensional analog of an ordinary
directed graph.

The category $\Lambda$ has range and source maps $r$ and $s$ (so
$\lambda$ is a morphism from $s(\lambda)$ to $r(\lambda)$).  For each
object $v$, and each $n \in \bbN^k$, let
$\Lambda^n(v) = \{\lambda \in \Lambda \mid d(\lambda) = n, r(\lambda)
= v\}$.  We assume throughout this paper that each $\Lambda^n(v)$ is a 
finite, non-empty set.  (This is usually expressed by saying that
$\Lambda$ is row finite and has no sources.)

A higher rank graph \cstar algebra, $C^*(\Lambda)$, is the universal
\cstar algebra generated by a family of partial isometries
$\{ s_{\lambda} \mid \lambda \in \Lambda \}$ satisfying:
\begin{enumerate}
\item $\{s_v \mid v \in \Lambda^0\}$ is a family of mutually
  orthogonal projections,
\item $s_{\lambda \mu} = s_{\lambda}s_{\mu}$, for all composable
$\lambda, \mu \in \Lambda$ (i.e., for all $\lambda,\mu$ with
$r(\mu)=s(\lambda)$),
\item $\sstars{\lambda}{\lambda} = s_v$ , where $v = s(\lambda)$,
\item for all $v \in \Lambda^0$ and all $n \in \bbN^k$,
$\ds s_v = \sum_{\lambda \in \Lambda^n(v)} \ssstar{\lambda}{\lambda}$.
\end{enumerate}
Any set of partial isometries in a \cstar algebra which satisfies
these four conditions is known as a Cuntz-Krieger family; if
$\{ t_{\lambda} \mid \lambda \in \Lambda\}$ is a Cuntz-Krieger family, 
then the map $s_{\lambda} \mapsto t_{\lambda}$ extends to a
homomorphism of $C^*(\Lambda)$ to the \cstar algebra generated by the
$t_{\lambda}$.

The description above is take largely from \cite{MR2001b:46102}, where 
the reader can find more detail and a number of examples.  The same
source provides more complete information about the path groupoid,
$\sG$, which we now summarize.

Let $\Omega_k$ denote the following $k$-graph:
\begin{itemize}
\item $\Obj \Omega_k = \bbN^k$.
\item $\Omega_k = \{(m,n) \mid (m,n) \in \bbN^k \times \bbN^k 
  \text{ and } m \leq n \}$.
\item $r(m,n)=m$; $s(m,n)=n$.
\item $d \colon \Omega_k \to \bbN^k$ by $d(m,n)=n-m$.
\end{itemize}

Infinite path space in $\Lambda$ is then defined to be
\begin{displaymath}
  \Lambda^{\infty} = \{ x \colon \Omega_k \to \Lambda \mid
x \text{ is a $k$-graph morphism}\}.
\end{displaymath}
For $v \in \Lambda^0$, let $\Lambda^{\infty}(v) =
  \{ x \in \Lambda^{\infty} \mid x(0) = v \}$.  For each 
$p \in \bbN^k$, define a shift map, 
$\sigma^p \colon \Lambda^{\infty} \to \Lambda^{\infty}$,  by 
$\sigma^p(x)(m,n) = x(m+p,n+p)$.

Using the factorization property, Kumjian and Pask show that
$x \in \Lambda^{\infty}$ is determined by the values $x(0,m)$, 
$m \in \bbN^k$.  
They also show that if
$\lambda \in \Lambda$ and $x \in \Lambda^{\infty}$ with 
$x(0) = s(\lambda)$, then
we can concatenate $\lambda$ and $x$: there is a
unique $y \in \Lambda^{\infty}$ such that 
$x = \sigma^{d(\lambda)} y$ and $\lambda = y(0,d(\lambda))$. 
 Naturally, we
write $y = \lambda x$.  This leads immediately to the factorization of 
any infinite path $x \in \Lambda^{\infty}$ as a product of a finite
path (an element of $\Lambda$) and an infinite tail:
$x = x(0,p)\sigma^p x$, for any $p \in \bbN^k$.

For any $\lambda \in \Lambda$, let
\begin{align*}
Z(\lambda) &= \{ \lambda x \in \Lambda^{\infty} \mid s(\lambda) = x(0) \} \\
&= \{ y \in \Lambda^{\infty} \mid y(0,d(\lambda)) = \lambda \}
\end{align*}
The collection $\{Z(\lambda) \mid \lambda \in \Lambda \}$ generates a
topology on path space $\Lambda^{\infty}$; in this topology each
$Z(\lambda)$ is a compact, open set.  The map $\lambda x \mapsto x$ is 
a homeomorphism of $Z(\lambda)$ onto $Z(s(\lambda))$ and each map
$\sigma^p$ is a local homeomorphism. 

$\Lambda^{\infty}$ will be identified with the set of
units in the groupoid $\sG_{\Lambda}$, which is defined by
\begin{displaymath}
  \sG_{\Lambda} = \{ (x,n,y) \in \Lambda^{\infty} \times \bbZ^k
\times \Lambda^{\infty} \mid \sigma^p x = \sigma^q y \text{ and }
n = p-q \}.
\end{displaymath}
When $k=1$, $\Lambda^{\infty}$ reduces to the usual infinite path
space and $\sG_{\Lambda}$ is the usual groupoid based on shift
equivalence on path space.  Inversion in $\sG_{\Lambda}$ is given by 
$(x,n,y)^{-1} = (y, -n, x)$.  Composable elements consist of those
with matching third and first coordinates, in which case multiplication 
is given by
$(x,n,y)(y,m,z) = (x,n+m,z)$.  
$\Lambda^{\infty}$ is identified
with the space of units, $\sG^0_{\Lambda}$, via
$x \mapsto (x,0,x)$.  A basis for a topology on $\sG_{\Lambda}$ is
given by the family
\begin{displaymath}
  Z(\lambda,\mu) = \{ (\lambda z, d(\lambda) - d(\mu), \mu z) \mid
z \in \Lambda^{\infty}(v) \},
\end{displaymath}
where $\lambda, \mu \in \Lambda$ and $s(\lambda) = s(\mu) = v$.
The topology generated by this basis is locally compact and
Haussdorff.  $\sG_{\Lambda}$ is then a second countable,
$r$-discrete, locally compact groupoid; each basic open set 
$Z(\lambda, \mu)$ is compact.  The identification of
$\Lambda^{\infty}$ with $\sG^0_{\Lambda}$ is a homeomorphism.  The
groupoid \cstar algebra, $C^*(\sG_{\Lambda})$, is isomorphic to the
higher rank graph \cstar algebra, $C^*(\Lambda)$.

The gauge action which appears in the spectral theorem for bimodules
is an action of the $k$-torus $\bbT^k$ on $C^*(\Lambda)$.  First, a
bit of notation: if $t \in \bbT^k$ and $n \in \bbN^k$ then 
$t^n = t_1^{n_1} t_2^{n_2} \cdots t_k^{n_k}$.  If 
$\{ s_{\lambda} \mid \lambda \in \Lambda \}$ is a generating
Cuntz-Krieger family, then so is
$\{ t^{d(\lambda)}s_{\lambda} \mid \lambda \in \Lambda \}$;  the
universal property then yields an automorphism $\gamma_t$ of
$C^*(\Lambda)$ such that 
$\gamma_t(s_{\lambda}) = t^{d(\lambda)}s_{\lambda}$, for all $\lambda$.

The fixed point algebra of the gauge action is an AF subalgebra of 
$C^*(\Lambda)$; it is generated by all $\ssstar{\lambda}{\mu}$ with
$d(\lambda)=d(\mu)$.  The map
$\Phi_0$ of $C^*(\Lambda)$ onto the fixed point algebra  given by
$\Phi_0(f) = \int_{\bbT^k} \gamma_t(f)\,dt$ is a faithful conditional
expectation.  For details concerning this, see \cite{MR2001b:46102}.

It is shown in \cite{MR2001b:46102} that $\sG_{\Lambda}$ is amenable;
consequently, $C^*(\sG_{\Lambda}) = C^*_{\text{red}}(\sG_{\Lambda})$.
Proposition II.4.2 in \cite{MR82h:46075} allows us to identify the
elements of $C^*(\sG_{\Lambda})$ with (some of the) elements of
$C_0(\sG_{\Lambda})$, the continuous functions on $\sG_{\Lambda}$
  vanishing at infinity.  (Note, however, that all continuous
  functions on $\sG_{\Lambda}$ with compact support are elements of
$C^*(\sG_{\Lambda})$.)

For each $m \in \bbZ^k$, let $\sG_m$ be the set of those elements
$(x,n,y)$ in $\sG_{\Lambda}$ with $n=m$.  The conditional expectation
$\Phi_0$ is just restriction map to $\sG_0$.  Restriction to $\sG_m$
is also a map of $C^*(\sG_{\Lambda})$ into itself; this is seen by
observing that it is given by the norm decreassing map $\Phi_m$
defined by 
$\Phi_m(f) = \int_{\bbT^k} t^{-m}\gamma_t(f)\,dt$.  If $\sB$ is a
closed linear subspace of $C^*(\Lambda)$ which is left invariant by
the gauge automorphisms, then $\Phi_m(\sB) \subseteq \sB$, for each $m$.

\section{The spectral theorem for bimodules}\label{S:stb}

Throughout this section, $\Lambda$ is a $k$-graph for which each
$\Lambda^n(v)$ is finite and non-empty and $\sG$ is the associated
$r$-discrete locally compact groupoid. Elements of the groupoid \cstar 
algebra (= higher rank graph \cstar algebra) are viewed as continuous
functions on $\sG$.
 (Since $k$ does not vary, 
we drop the subscript from the notation for the grouoid.)  As above,
we identify path space $\Lambda^{\infty}$ with the space of units of
$\sG$; with this identification  $C_0(\Lambda^{\infty})$ becomes an
abelian subalgebra of $C^*(\sG)$.  $\Lambda^{\infty}$ is not compact
except when $\Lambda$ has finitely many objects (``vertices''), hence
the use of $C_0$.

For simplicity of notation, let $\sA$ denote the groupoid \cstar
algebra and let $\sD$ denote $C_0(\Lambda^{\infty})$.
At the end of the section we will discuss when $\sD$ is a masa in
$\sA$.

Since $\sG$ is $r$-discrete, the Haar system can be taken to be
counting measure, and so is not mentioned explicitly.  Since elements
of $\sA$ are interpreted as functions on $\sG$, multiplication is
given by a convolution type formula
\begin{displaymath}
  fg(x,n,y) = \sum f(x,p,z)g(z,q,y)
\end{displaymath}
where the sum is taken over all composable pairs $(x,p,z)$ and
$(z,q,y)$ with $p+q=n$.  (For functions in $\sA$, the series will
converge.)  In particular, if $f \in \sA$ and $g \in \sD$,
\begin{align} \label{ADproduct1}
gf(x,n,y) &= g(x,0,x)f(x,n,y), \\
gf(x,n,y) &= f(x,n,y)g(y,0,y). \label{ADproduct2}
\end{align}

For each $\lambda \in \Lambda$, let $s_{\lambda}$ denote the
characteristic function of the set $Z(\lambda,s(\lambda))$.  Then
$\{ s_{\lambda} \mid \lambda \in \Lambda \}$ forms a Cuntz-Krieger
family and generates $\sA$ as a \cstar algebra.  This can be 
checked using the definition of $Z(\lambda,s(\lambda))$ and the
formula given above for multiplication.  Note also that, for 
$\lambda, \mu \in \Lambda$ with $s(\lambda)=s(\mu)$, 
$\ssstar{\lambda}{\mu}$ is the characteristic function of
the set $Z(\lambda,\mu)$.

If $\sB \subseteq \sA$ is a bimodule over $\sD$, we define the
\emph{spectrum} of $\sB$ to be:
\begin{displaymath}
  \sigma(\sB) = \{(x,n,y) \in \sG \mid f(x,n,y) \ne 0 \text{ for some } 
  f \in \sB\}.
\end{displaymath}
The spectrum $\sigma(\sB)$ is an open subset of $\sG$.  On the other
hand, any open subset $P$ of $\sG$ determines a $\sD$-module $A(P)$
given by
\begin{displaymath}
  A(P) = \{f \in \sA \mid f(x,n,y) =0 \text{ for all } (x,n,y) \notin
  P \}
\end{displaymath}

Since $P$ is open, if $(x,n,y) \in P$, then there is a basic open set
$Z(\lambda, \mu)$ such that 
$(x,n,y) \in Z(\lambda, \mu) \subseteq P$.  It follows that
$\ssstar{\lambda}{\mu} \in A(P)$; since $\ssstar{\lambda}{\mu}$ has the
value 1 at $(x,n,y)$, we obtain $\sigma(A(P)) =P$, for any open subset 
$P \subseteq \sG$.

It is clear that if $\sB$ is a bimodule over $\sD$ then
$\sB \subseteq A(\sigma(\sB))$; equality does not always hold.  A
counterexample in the special case of Cuntz algebras (algebras
determined by 1-graphs with only one vertex) can be found in
\cite{hop_peters}.  Also, it is shown in \cite{hpp} that there is a
counterexample for any graph \cstar algebra which is not AF.  (For AF
\cstar algebras the Muhly-Solel spectral theorem for bimodules says
that $\sB = A(\sigma(\sB))$ always.)  Thus counterexamples exist for
all 1-graphs which contain a loop.

A characterization of those bimodules which are determined by their
spectra -- $\sB =A(\sigma(\sB))$ -- is given in the graph \cstar
algebra context in \cite{hpp}.  The main result in this note is the
extension  to the higher rank context:

\begin{theorem}[Spectral Theorem for Bimodules]
Let $\Lambda$ be a row finite $k$-graph with no sources.  Let
$\sG$ be the associated path groupoid.
Let $\sA = C^*(\Lambda) = C^*(\sG)$ and 
$\sD = C_0(\Lambda^{\infty})$.  If $\sB \subseteq \sA$ is a bimodule
over $\sD$, then the following are equivalent:
\begin{enumerate}
\item \label{stb_det}  $\sB = A(\sigma(\sB))$.
\item \label{stb_ck} $\sB$ is generated by the Cuntz-Krieger partial
  isometries which it contains.
\item \label{stb_gauge}  $\sB$ is invariant under the gauge automorphisms.
\end{enumerate}
\end{theorem}

\begin{proof}
$(\ref{stb_det}) \Rightarrow (\ref{stb_ck})$.  Assume $P$ is an open
subset of $\sG$.  Let $\sB$ be the bimodule generated by the
Cuntz-Krieger partial isometries in $A(P)$.  Each such partial
isometry has its support in $P$, so $\sigma(\sB) \subseteq P$ and 
$\sB \subseteq A(P)$.  We need to show that any function $f$ in
$A(P)$ is actually in $\sB$.  We claim that it is sufficient to do
this for functions which are supported on some
$Z(\lambda,\mu) \subseteq P$.  Indeed, it then follows readily that
functions supported on compact subsets of $P$ are in $\sB$ (every
compact subset of $P$ is contained in a finite union of subsets of the 
form $Z(\lambda,\mu)$) and the compactly supported functions in $A(P)$ 
are dense in $A(P)$.

If $f$ has support in $Z(\lambda,\mu)$, with the aid of convolution formulas
(\ref{ADproduct1}) and (\ref{ADproduct2})
 it is easy to find a function $g$ supported in
$\Lambda^{\infty}$ such that $f = g \ssstar{\lambda}{\mu}$.  Since
$\ssstar{\lambda}{\mu} \in \sB$ and $g \in \sD$, $f \in \sB$ also.

$(\ref{stb_ck}) \Rightarrow (\ref{stb_gauge})$.  Since a gauge
automorphism maps a Cuntz-Krieger partial isometry to a scalar
multiple of itself, $\sB$ is trivially left invariant when it is
generated by its Cuntz-Krieger partial isometries.

$(\ref{stb_gauge}) \Rightarrow (\ref{stb_det})$.  Let $\sB$ be a gauge 
invariant bimodule and let $P = \sigma(\sB)$.  Since 
$\sB \subseteq A(P)$ is automatic, we just need to show that
$A(P) \subseteq \sB$.  For each $m \in \bbZ^k$, let
$P_m = P \cap \sG_m$, so that $P = \cup_m P_m$.  Since
$\Phi_m$ maps $A(P)$ onto $A(P_m)$ and, for each $f$,
$f$ is in the closed linear span of the $\Phi_m(f)$,  we need merely show
that  $A(P_m) \subseteq \sB$, for each $m$.

Fix $m$.  Suppose that $\alpha, \beta \in \Lambda$ satisfy
$s(\alpha)=s(\beta)$ and $d(\alpha) - d(\beta) = m$.  Denote
$\sG_{\alpha,\beta} = \{(\alpha z, m, \beta w) \mid z,w \in
\Lambda^{\infty}(s(\alpha)) \}$ and
$P_{\alpha,\beta} = P_m \cap \sG_{\alpha,\beta}$.  Now, by  
what we have just proven $A(P_m)$ is the 
closed linear span of the Cuntz-Krieger partial isometries which it
contains.  But if $\ssstar{\alpha}{\beta}$ is one of these, then
$\ssstar{\alpha}{\beta} \in A(P_{\alpha,\beta})$, so $A(P_m)$ is the
closed linear span of the $A(P_{\alpha,\beta})$.  This reduces the
task to showing that $A(P_{\alpha,\beta}) \subset \sB$ for each
suitable pair $\alpha, \beta$.

We can finish the proof by transfering the problem to (a subset of)
$\sG_0$; the latter is a principal groupoid so the Muhly-Solel
spectral theorem for bimodules is available.  Let
\begin{displaymath}
  \sG_0(s(\alpha)) = \{(z,0,w) \mid z,w \in \Lambda^{\infty}(s(\alpha))\}.
\end{displaymath}
The map $\psi \colon \sG_0(s(\alpha)) \to \sG_{\alpha,\beta}$ given by
$(z,0,w) \mapsto (\alpha z, m, \beta w)$ is a homeomorphism.  Let
$Q$ be the inverse image of $P_{\alpha,\beta}$ under this map.  Note
that $f \mapsto s_{\alpha}^{\phantom{*}} f s^*_{\beta}$ carries
$A(Q)$ onto $A(P_{\alpha,\beta})$.

Let
\begin{displaymath}
  \sC = \{f \in A(\sG_0(s(\alpha))) \mid  s_{\alpha}^{\phantom{*}} f
  s^*_{\beta} \in \sB \}.
\end{displaymath}
We claim that $\sC$ is a bimodule over $\sD$.  Since $\sD$ is
generated by projections of the form $\ssstar{\lambda}{\lambda}$, it
suffices to show that $\sC$ is closed under multiplication left and
right by such projections.  Now if $f \in \sC$, then, since
$ s^{\phantom{*}}_{\alpha} \ssstar{\lambda}{\lambda} \ne 0$ exactly
when $\ssstar{\lambda}{\lambda} \leq \sstars{\alpha}{\alpha}$,
\begin{displaymath}
 s^{\phantom{*}}_{\alpha}  \ssstar{\lambda}{\lambda} f
 s^*_{\beta} =  s^{\phantom{*}}_{\alpha}
 \ssstar{\lambda}{\lambda} \sstars{\alpha}{\alpha} f
 s^*_{\beta} = \ssstar{\alpha \lambda}{\alpha \lambda} 
 s^{\phantom{*}}_{\alpha} f  s^*_{\beta} \in \sB.
\end{displaymath}
The last assertion uses 
$\ssstar{\alpha \lambda}{\alpha \lambda} \in \sD$ and
$ s^{\phantom{*}}_{\alpha} f  s^*_{\beta} \in \sB$.
Thus $\sC$ is a left bimodule over $\sD$; the argument that it is a
right bimodule is similar.

The definition of $Q$ implies that $\sigma(\sC) \subseteq Q$.  The
gauge invariance of $\sB$ implies that $\sigma(\sC) = Q$.  Indeed, let
$q \in Q$ and let $p = \psi(q)$.  Since $p \in P$, there is $f \in
\sB$ such that $f(p) \neq 0$.  Then $\Phi_m(f)(p) \neq 0$ and,
by gauge invariance, $\Phi_m(f) \in \sB$.  If
$g = s^*_{\alpha} \Phi_m(f) s^{\phantom{*}}_{\beta}$, then
$g \in \sC$ and $g(q) \neq 0$.

Since $\sigma(\sC) = Q$ and the Muhly-Solel spectrum for bimodules
holds in $A(\sG_0)$, we have $\sC = A(Q)$.  This implies that
$A(P_{\alpha,\beta}) \subseteq \sB$.
\end{proof}

As mentioned earlier, $\sD = C_0(\Lambda^{\infty})$ need not be a masa
in $\sA$.  For the graph C*-algebra case, it was shown in \cite{hpp} 
that $\sD$ is a masa if, and only if, every loop has an entrance.  
Kumjian and Pask \cite{MR2001b:46102} define an analogous condition, the 
\emph{aperiodicity condition}, for higher rank graphs and use this to
extend the 
Cuntz-Krieger uniqueness theorem.  Their condition also extends the 
masa theorem.  Here are the relevant definitions: 
an element $x \in \Lambda^{\infty}$ is \emph{periodic} with non-zero period 
$p \in \bbZ^k$ if, for every $(m,n) \in \Omega$ with $m+p \geq 0$,
$x(m+p,n+p) = x(m,n)$.  If there is an element $n \in \bbN^k$ such
that $\sigma^n(x)$ is periodic, $x$ is \emph{eventually periodic};
otherwise, $x$ is \emph{aperiodic}.  Finally, $\Lambda$ satisfies the 
\emph{aperiodicity condition} if, for every $v \in \Lambda^0$, there
is an aperiodic path $x \in \Lambda^{\infty}(v)$.

Note that $x$ is eventually periodic with period $p$ if, and only if, 
$(x,p,x) \in \sG$.

Kumjian and Pask prove that $\Lambda$ satisfies the aperiodicity
condition if, and only if, the points in $\sG$ with trivial isotropy
are dense in $\sG^0$ \cite[Proposition 4.5]{MR2001b:46102}.
We will show below that the aperiodicity condition is also equivalent
to the assertion that $\sG^0$ is the interior of the isotropy group
bundle $\sG^1$.  (Note: in the Kumjian-Pask proposition, $\sG^0$ is
viewed as $\Lambda^{\infty}$; we will view $\sG^0$ as the open subset
$\{(x,0,x) \mid x \in \Lambda^{\infty}\}$ of
$\sG^1 = \{(x,p,x) \in \sG \mid p \in \bbZ^k\}$.)
Renault \cite[Proposition II.4.7]{MR82h:46075} 
has shown that, $C_0(\sG^0)$ is a masa in 
$C^*_{\text{red}}(\sG)$ if, and only if, $\sG^0$ is the interior of
$\sG^1$. Since the path groupoid $\sG$ is amenable, Renault's
Proposition yields the masa theorem.

\begin{proposition}
$\Lambda$ satisifies the aperiodicity condition if, and only if,
$\sG^0$ is the interior of $\sG^1$.
\end{proposition}

\begin{proof}
Assume that the aperiodicity condition holds.  Let $(x,p,x) \in \sG^1$
with $p \neq 0$.  We shall show that we can approximate $(x,p,x)$ by
points in $\sG$ which are not in $\sG^1$.  This shows that $(x,p,x)$
is not in the interior of $\sG^1$.  Since $\sG^0$ is an open subset of
$\sG^1$, it follows that $\sG^0$ is the interior.

Let $Z(\alpha,\beta)$ be a neighborhood of $(x,p,x)$.  For $m$
sufficiently large (meaning for each $m_i$ sufficiently large), 
$m+p \geq 0$ and both $x(0,m)$ and $x(0,m+p)$ lie in $Z(\alpha)$ and
in $Z(\beta)$.  Since $(x,p,x) \in \sG$, $\sigma^m(x) =
\sigma^{m+p}(x)$ and $x(0,m)$ and $x(0,m+p)$ have a
common source $v$.  Choose $y$ aperiodic
in $\Lambda^{\infty}(v)$.  Let $z = x(0,m)y$ and $w = x(0, m+p)y$.
Then $z \neq w$ and $(z,p,w) \in Z(\alpha,\beta)$.  So 
$(z,p,w) \notin \sG^1$ and $(z,p,w)$ approximates $(x,p,x)$.

Now suppose that $\Lambda$ does not satisfy the aperiodicity
condition.  By Proposition 4.5 in \cite{MR2001b:46102}, there is 
$x \in \Lambda^{\infty}$ which cannot be approximated by aperiodic
points.  Since $x$ must be  eventually periodic
there is a non-zero element $p$ 
of $\bbZ^k$ such that $(x,p,x) \in \sG$.  If $(x,p,x)$ could be
approximated in the topology of $\sG$ by points outside $\sG^1$, it
would follow that $x$ is a limit of aperiodic points in
$\Lambda^{\infty}$ -- a contradiction.  This shows that $(x,p,x)$ is
in the interior of $\sG^1$ and so $\sG^0$ is not the interior.
\end{proof}

This Proposition, Proposition II.4.7 in \cite{MR82h:46075}, and the
amenability of $\sG$ yield the following theorem.

\begin{theorem}
$\sD$ is a masa in $\sA$ if, and only if $\Lambda$ satisifies the
aperiodicity condition.  
\end{theorem}

\providecommand{\bysame}{\leavevmode\hbox to3em{\hrulefill}\thinspace}
\providecommand{\MR}{\relax\ifhmode\unskip\space\fi MR }
\providecommand{\MRhref}[2]{%
  \href{http://www.ams.org/mathscinet-getitem?mr=#1}{#2}
}
\providecommand{\href}[2]{#2}

\end{document}